\documentclass[a4paper,twoside,10pt,reqno]{article}
\usepackage[english]{babel}
\usepackage[dvips]{graphicx}
\usepackage[T1]{fontenc}
\usepackage[latin1]{inputenc}
\usepackage{amssymb,amsthm,dsfont,mathrsfs}
\usepackage{amsfonts,amsmath,euscript}
\usepackage{color}

\newtheorem{teo}{Theorem}[section]
\newtheorem{pro}[teo]{Proposition}

\newtheorem{rem}{Remark}[section]

\newcounter{example}[section]

           \addto{\captionsenglish}{}

\newcommand{\hs}{\hspace{3pt}}
\newcommand{\dst}{\displaystyle}
\newcommand{\dem}{{\bf Dem. }}
\newcommand{\fdem}{$\square$}

\newcommand{\sscr}{\scriptscriptstyle}

\newcommand{\nn}{\nonumber}
\newcommand{\nid}{\noindent}

\newcommand{\titulo}[1]{\mbox{} \\ \noindent \textit{\textbf{\Large #1}}\\}

\renewcommand{\abstract}[1]{{\small \noindent \textbf{Abstract:} #1\\}}

\newcommand{\keywords}[1]{{\small \noindent \textbf{Keywords:} #1\\}}

\begin{document}

\titulo{On the penultimate tail behavior of Weibull-type models}

%\begin{center}
%\autor{Marta Ferreira}
%
%\afil{Department of Mathematics, University of Minho, Portugal\\
%msferreira@math.uminho.pt\\}
%\end{center}

 \abstract{The Gumbel max-domain of attraction corresponds to
a null tail index which do not distinguish the different tail
weights that might exist between distributions within this class.
The Weibull-type distributions form an important subgroup of this
latter and includes the so-called \emph{Weibull-tail coefficient},
usually denoted $\theta$, that specifies the tail behavior, with
larger values indicating slower tail decay. Here we shall see that
the Weibull-type distributions present a penultimate tail behavior
Fréchet if $\theta>1$ and a penultimate tail behavior Weibull
whenever $\theta<1$. }


%\nid\textbf{AMS 2000 Subject Classification} 60G70\\

\section{Introduction}\label{sint}
The main objective of an extreme value analysis is to estimate the
probability of events that are more extreme than any that have
already been observed. By way of example, suppose that a sea-wall
projection requires a coastal defense from all sea-levels, for the
next $100$ years. The use of extremal models enables extrapolations
of this type. The central result in classical Extreme Value Theory
(EVT) states that, for an i.i.d.\hs sequence, $\{X_{n}\}_{\sscr
n\geq 1}$, having common distribution function (d.f.) $F$, if there
are real constants $a_{n}>0$ and $b_{n}$ such that,
\begin{eqnarray}\label{domatrac}
P(\max(X_{1},...,X_{n})\leq
a_{n}x+b_{n})=F^n(a_{n}x+b_{n})\longrightarrow_{n\rightarrow\infty} G_{\gamma}(x)\, ,
\end{eqnarray}
for some non degenerate function $G_{\gamma}$, then it must be the
Generalized Extreme Value function (\textit{GEV}),
\begin{eqnarray}\label{gev}\nn
G_{\gamma}(x)= \exp(-(1+\gamma x)^{-1/\gamma})\textrm{, }1+\gamma
x>0\textrm{, }\gamma \in \mathbb{R},
\end{eqnarray}
($G_0(x)=\exp(-e^{-x})$) and we say that $F$ belongs to the
max-domain of attraction of $G_\gamma$, in short, $F\in
\mathcal{D}(G_{\gamma})$. The parameter $\gamma$, known as the tail
index, is a shape parameter as it determines the tail behavior of
$F$, being so a crucial issue in EVT. More precisely, if $\gamma>0$
we are in the domain of attraction Fréchet corresponding to a heavy
tail, $\gamma<0$ indicates the Weibull domain of attraction of light
tails and $\gamma=0$ means a Gumbel domain of attraction and an
exponential tail.

However, as Fisher and Tippett (\cite{fish+tip}, 1928) remarked, if
one approximates the distribution of the successive maxima of normal
samples not by the limit distribution Gumbel but by a sequence of
other extreme value distributions converging to the limit
distribution, the approximation is asymptotically improved. They
called penultimate distributions to this sequence of approximating
extreme value distributions.

Here we will analyze the penultimate tail behavior of the
Weibull-type models which are a wide class in the Gumbel
max-domain.\\

 Weibull-type models have representation
\begin{eqnarray}\label{wei}
1-F(x)=\exp(-H(x)),\,x\geq x_0\geq 0,\,\theta>0%\textrm{, where $H(x)=x^{1/\theta}l(x)$ or $H^{-1}(x)=x^{\theta}l^*(x)$,}
\end{eqnarray}
where
\begin{eqnarray}\label{H}
\textrm{$H(x)=x^{1/\theta}l(x)$ or
$H^{-1}(x)=x^{\theta}l^*(x)$,}
\end{eqnarray}
with $H^{-1}$ denoting the generalized inverse of $H$ and functions
$l$ and $l^*$ are slowly varying at infinity (i.e., $l(tx)/l(t)\to
1$ as $t\to\infty$ for all $x>0$ and the same holds for $l^*$).

We say that functions $H$ and $H^{-1}$ are regularly varying with
indexes, $1/\theta$ and $\theta$, respectively.

 The parameter
$\theta$, called the Weibull-tail coefficient, governs the tail
behavior of $F$, with larger values indicating slower tail decay.
The Weibull-type distributions form an important subgroup within the
Gumbel class ($\gamma=0$) an the tail behavior can then be specified
using the Weibull-tail coefficient (Dierckx \emph{et al.},
\cite{dierckx+} 2009).

Weibull-type models include well-known distributions such as, Normal
($\theta=1/2$), Weibull($\alpha$,$\lambda$) ($\theta=1/\alpha$),
Extended Weibull($\beta$,$\delta$) ($\theta=1/\beta$), Exponential,
Gamma and Logistic
($\theta=1$) (Gardes and Girard \cite{gardes+girard}, 2008).\\

%Here we shall analyze the penultimate tail behavior of Weibull-type
%models in (\ref{wei}) through the class in (\ref{weilcc}), since the
%conditions of the results that we are going to use only concern
%limits as $x$ goes to infinity.

We will prove that the Weibull-type models present penultimate tail
behavior Fréchet or Weibull whenever $\theta>1$ or $\theta<1$,
respectively. To this end, we will use a result in de Haan and Gomes
(\cite{gomes+haan}, 1999) and in Gomes (\cite{gomes}, 1984).\\

We remark that the case $\theta=1$ cannot be treated in a unified
way since we can find penultimate tail indexes with different
orders. For instance, the Exponential and Logistic have $\gamma_n$
of order $1/n$ but the Gamma distribution have $\gamma_n$ of order
$1/\log^2 n$ (Gomes \cite{gomes2}, 1993). Therefore, we shall assume
that $\theta\not=1$ all over the paper.

\section{Results}

\nid Consider
\begin{eqnarray}\label{defk}
k(x):=\frac{d}{dx}[-\log(-\log F(x))]%=-\frac{F'(x)}{F(x)\log F(x)}=H'(x).
\end{eqnarray}
%Observe that
%\begin{eqnarray}\label{weilcc}
%F^{-1}(x)=H^{-1}(-\log(-\log x))=(-\log(-\log x))^{\theta}l^*(-\log(-\log x))%\textrm{, where $H(x)=x^{1/\theta}l(x)$ or $H^{-1}(x)=x^{\theta}l^*(x)$}
%\end{eqnarray}
%

Observe that
$$
-\log F(x)=(1-F(x))(1+O(1-F(x)))
$$
and, by (\ref{wei}),
$$
-\log(-\log F(x))=H(x)-\log (1+O(e^{-H(x)})),
$$
leading to
\begin{eqnarray}\label{k}
k(x)=H'(x)\Big(1+\frac{O(e^{-H(x)})}{1+O(e^{-H(x)})}\Big)=H'(x)(1+o(1)).
\end{eqnarray}
Therefore
\begin{eqnarray}\label{kl}
\begin{array}{rl}
k'(x)=&H''(x)\Big(1+\frac{O(e^{-H(x)})}{1+O(e^{-H(x)})}\Big)
+[H'(x)]^2\frac{O(e^{-2H(x)})-O(e^{-H(x)})(1+O(e^{-H(x)}))}{(1+O(e^{-H(x)}))^2}\vspace{0.35cm}\\
=&H''(x)(1+o(1)).
\end{array}
\end{eqnarray}

Analogously we derive
\begin{eqnarray}\label{kll+lll}
k''(x)=H'''(x)(1+o(1))\textrm{ and }k'''(x)=H^{(iv)}(x)(1+o(1)).
\end{eqnarray}

%Thus being, we can apply (\ref{ks}) to the Weibull-type model in
%(\ref{wei}).

Consider $H(x)$ given in (\ref{H}). We have
\begin{eqnarray}%\label{ks}
\begin{array}{rl}
H'(x)=&  x^{\theta^{-1}-1}l(x)\Big(\theta^{-1} +\frac{xl'(x)}{l(x)}\Big) \vspace{0.35cm}\\
 H''(x)=&  x^{\theta^{-1}-2}l(x)\Big(\theta^{-1}(\theta^{-1}-1) +2\theta^{-1}\frac{xl'(x)}{l(x)}+\frac{x^{2}l''(x)}{l(x)}\Big) \vspace{0.35cm}\\
 H'''(x)=&  x^{\theta^{-1}-3}l(x)\Big(\theta^{-1}(\theta^{-1}-1)(\theta^{-1}-2)
 +3\theta^{-1}(\theta^{-1}-1)\frac{xl'(x)}{l(x)}+3\theta^{-1}\frac{x^{2}l''(x)}{l(x)}
 +\frac{x^{3}l'''(x)}{l(x)}\Big) \vspace{0.35cm}\\
 H^{(iv)}(x)=&  x^{\theta^{-1}-4}l(x)\Big(\theta^{-1}(\theta^{-1}-1)(\theta^{-1}-2)(\theta^{-1}-3)
 +4\theta^{-1}(\theta^{-1}-1)(\theta^{-1}-2)\frac{xl'(x)}{l(x)}\vspace{0.35cm}\\
 &\hspace{1.7cm}+4\theta^{-1}(\theta^{-1}-1)\frac{x^{2}l''(x)}{l(x)}
 +\theta^{-1}\frac{x^{3}l'''(x)}{l(x)}+\frac{x^{4}l^{(iv)}(x)}{l(x)}\Big).
\end{array}
\end{eqnarray}
Assuming that, as $x\to\infty$,
\begin{eqnarray}\label{condls}
\begin{array}{cccc}
\frac{xl'(x)}{l(x)}\to 0,  &\frac{x^2l''(x)}{l(x)}\to 0,&\frac{x^3l'''(x)}{l(x)}\to 0\,\,\textrm{ and}&\frac{x^4l^{(iv)}(x)}{l(x)}\to 0,
\end{array}
\end{eqnarray}
and applying (\ref{k})-(\ref{kll+lll}), we have that
\begin{eqnarray}\label{ks}
\begin{array}{l}
x k(x)=   \theta^{-1} x^{\theta^{-1}}l(x)(1+o(1))\\
x^2k'(x)=x k(x)(\theta^{-1}-1)\\%= \theta^{-1} x^{\theta^{-1}}l(x)(\theta^{-1}-1)(1+o(1))\\
x^3k''(x)=x k(x)(\theta^{-1}-2)(\theta^{-1}-1)\\%= \theta^{-1} x^{\theta^{-1}}l(x)(\theta^{-1}-2)(\theta^{-1}-1)(1+o(1))\\
  x^4k'''(x)=x k(x)(\theta^{-1}-3)(\theta^{-1}-2)(\theta^{-1}-1)\\%= \theta^{-1} x^{\theta^{-1}}l(x)(\theta^{-1}-3)(\theta^{-1}-2)(\theta^{-1}-1)(1+o(1)).\\
\end{array}
\end{eqnarray}
Note that, if $l(x)$ is monotone for $x\geq x_0>0$, then
$xl'(x)/l(x)\to 0$, as $x\to\infty$. The other conditions in
(\ref{condls}) are also satisfied by the most common models.\\

Observe that equations (\ref{k})-(\ref{kll+lll}) hold for the class
of d.f.'s
\begin{eqnarray}\label{weilcc}
-\log F(x)=\exp(-H(x))%\textrm{, where $H(x)=x^{1/\theta}l(x)$ or $H^{-1}(x)=x^{\theta}l^*(x)$}
\end{eqnarray}
which have been studied in Canto e Castro (\cite{lcc}, 1992). More
precisely, it was derived the structure of the remainder
$F^n(a_nx+b_n)-G_{\gamma}(x)$, by showing that $F$ is in class
$\mathbf{A_1}$ of Anderson (\cite{anderson}, 1976), i.e.,
\begin{eqnarray}\label{p344lcc}
\lim_{x\to\infty}\frac{k''(x)}{k(x)k'(x)}=\lim_{x\to\infty}\frac{x^3k''(x)}{xk(x)x^2k'(x)}=\lim_{x\to\infty}\frac{(\theta^{-1}-2)}{xk(x)}=0.
\end{eqnarray}
leading to
\begin{eqnarray}\label{rateconv}
F^n(a_nx+b_n)-G_{\gamma}(x)=\frac{x^2}{2}\Big(\frac{k'(b_n)}{k^2(b_n)}\Big)g_{\gamma}(x)(1+o(1)),
\end{eqnarray}
uniformly for $x$ in bounded intervals in the support of
$G_{\gamma}$, where $g_{\gamma}(x)=G'_{\gamma}(x)$, $a_n=1/k(b_n)$
and $F(b_n)=\exp(-1/n)$. More details can be seen in Canto e Castro
(\cite{lcc}, 1992).\\
%Proposition 3.4.4).\\

Observe that
\begin{eqnarray}\label{bn}
b_n=H^{-1}(\log n),%=(\log n)^{\theta}l^*(\log n).
\end{eqnarray}
and the rate of convergence of $F^n(a_nx+b_n)$ to $G_{\gamma}(x)$ is
\begin{eqnarray}\label{kbn}
\frac{k'(b_n)}{k^2(b_n)}= \frac{b_n^2k'(b_n)}{(b_nk(b_n))^2} \sim
\frac{\theta^{-1}-1}{\theta^{-1}(\log n)}=\frac{1-\theta}{\log n}.
 \end{eqnarray}
%For larger values of $x$, $1-F(x)$ is asymptotically equivalent to
%$-\log F(x)$, and hence, the tail behavior of $F$ can be analyzed
%through the latter.

\bigskip

\bigskip

Now consider $x^F:=\sup\{x:F(x)<1\}$ and
\begin{eqnarray}\label{phi}
\varphi(t)=(1/k)'(t)=-k'(t)/(k(t))^2.
\end{eqnarray}
Based on the von Mises' first order condition,
\begin{eqnarray}\label{1vmis}
\dst\lim_{t\to x^F}\varphi(t)=0,
\end{eqnarray}
the von Mises' second order condition,
\begin{eqnarray}\label{2vmis}
\dst\lim_{t\to x^F}\frac{\varphi'(t)}{k(t)\varphi(t)}=0,
\end{eqnarray}
and von Mises' type penultimate condition,
\begin{eqnarray}\label{penvmis}
\dst\lim_{t\to x^F}\frac{\varphi''(t)}{k(t)\varphi'(t)}=0,
\end{eqnarray}
Theorem 1 in Haan and Gomes (\cite{gomes+haan}, 1999) allow us to
derive bounds for $(F^n(a_n x+b_n)-G_{\gamma_n}(x))/\gamma'(\log
n)$, with $a_n$ and $b_n$ given above, where
\begin{eqnarray}\label{gamat}
\dst\gamma(t)=\varphi(H^{-1}(t))=-\frac{k'(H^{-1}(t))}{k^2(H^{-1}(t))}
\end{eqnarray}
and the penultimate tail index is
\begin{eqnarray}\label{gaman}
\gamma_n=\gamma(\log n)=-k'(b_n)/k^2(b_n).
\end{eqnarray}

%In the following we see that the von Mises' conditions
%(\ref{1vmis})-(\ref{penvmis}) are fulfilled by model in
%(\ref{weilcc}) (where $x^F=\infty$), by using equations in
%(\ref{ks}).

\begin{pro}\label{phaan+gomes}
Weibull-type models satisfy  the von Mises' conditions
(\ref{1vmis})-(\ref{penvmis}) whenever conditions in (\ref{condls})
are fulfilled. Moreover, they present penultimate tail behavior
Fréchet if $\theta>1$ and Weibull if $\theta<1$.
\end{pro}
\dem First observe that $x^F=+\infty$,
$$
\begin{array}{c}
\dst\varphi'(t)=-\frac{k''(t)k(t)-2(k'(t))^2}{(k(t))^3}
\end{array}
$$
and
$$
\begin{array}{c}
\dst\varphi''(t)=6\frac{k'(t)}{(k(t))^2}\Big(\frac{k''(t)}{k(t)}-\Big(\frac{k'(t)}{k(t)}\Big)^2\Big)-\frac{k'''(t)}{(k(t))^2}.
\end{array}
$$
Now we have, successively,
\begin{eqnarray}\label{1vmisres}
\dst\lim_{t\to x^F}\varphi(t)=\dst\lim_{t\to x^F}-\frac{k'(t)}{(k(t))^2}=\dst\lim_{t\to x^F}-\frac{t^2k'(t)}{(tk(t))^2}
=\dst\lim_{t\to x^F}-\frac{(\theta^{-1}-1)(1+o(1))}{\theta^{-1} t^{\theta^{-1}}l(t)(1+o(1))}=0,
\end{eqnarray}
by (\ref{p344lcc}) and (\ref{1vmisres}),
\begin{eqnarray}\label{2vmisres} \dst\lim_{t\to
x^F}\frac{\varphi'(t)}{k(t)\varphi(t)}=\lim_{t\to
x^F}\Big(2\frac{k'(t)}{(k(t))^2}-\frac{k''(t)}{k(t)k'(t)}\Big)=0,
\end{eqnarray}
and, applying (\ref{ks}),
\begin{eqnarray}\label{penvmisres}
\begin{array}{rl}
&\dst\lim_{t\to x^F}\frac{\varphi''(t)}{k(t)\varphi'(t)}
=\dst\lim_{t\to x^F}\frac{6k'(t)k''(t)k(t)-6(k'(t))^3-k'''(t)(k(t))^2}{(k(t))^2(2(k'(t))^2-k''(t)k(t))} \vspace{0.35cm}\\
=&\dst\lim_{t\to x^F}\frac{6(\theta^{-1}-1)(\theta^{-1}-2)-6(\theta^{-1}-1)^2-(\theta^{-1}-2)(\theta^{-1}-3)}{tk(t)[2(\theta^{-1}-1)-(\theta^{-1}-2)]}=0.
\end{array}
\end{eqnarray}
 By (\ref{kbn}), the penultimate tail index, $\gamma_n$, in (\ref{gaman}) is given by
\begin{eqnarray}\label{gamanc}
\gamma_n=\gamma(\log n)=-k'(b_n)/k^2(b_n)\sim\frac{\theta-1}{\log n}.
\end{eqnarray}
Therefore, we obtain $\gamma_n>0$ and $\gamma_n<0$ if, respectively,
$\theta>1$ and $\theta<1$.

\nid Now we compute the rate of convergence $\gamma'(\log n)$.
Observe that, after some calculations,
\begin{eqnarray}\label{gamatl}
\begin{array}{rl}
\dst\gamma'(t)=\frac{2k'(H^{-1}(t))^2-k''(H^{-1}(t))k(H^{-1}(t))}{k(H^{-1}(t))^4},
\end{array}
\end{eqnarray}
and applying (\ref{bn}) we have
\begin{eqnarray}\label{gamatl}
\begin{array}{rl}
\dst\gamma'(\log n)=\frac{2(b_nk(b_n)(\theta^{-1}-1))^2-(b_nk(b_n))^2(\theta^{-1}-1)(\theta^{-1}-2)}
{(b_nk(b_n))^4}\sim \frac{2\theta(1-\theta)}{(\log n)^2}.  \,\,$\fdem$
\end{array}
\end{eqnarray}
%Observe that the rate of convergence vanishes at $\theta=1$.\\

%\bigskip

\begin{rem}
The result in Proposition \ref{phaan+gomes} can also be derived
based on Gomes (\cite{gomes}, 1984). This latter imposes the
condition
\begin{eqnarray}\label{gomes84}
\begin{array}{c}
\dst\lim_{t\to x^F}\frac{\varphi'(t)}{k(t)(\varphi(t))^2}=c<\infty
\end{array}
\end{eqnarray}
where $\varphi(t)=(1/k)'(t)=-k'(t)/(k(t))^2$, and obtains
boundedness of $(F^n(a_n x+b_n)-G_{\gamma_n}(x))/\gamma_n^2$, with
penultimate tail index, $\gamma_n$, given in (\ref{gamanc}).

Observe that
$$
\begin{array}{c}
\dst\frac{\varphi'(t)}{k(t)(\varphi(t))^2}=2-\frac{k''(t)k(t)}{(k'(t))^2}.
\end{array}
$$
Applying (\ref{ks}), we obtain
$$
\begin{array}{c}
\dst\lim_{t\to x^F}\frac{k''(t)k(t)}{(k'(t))^2}=\lim_{t\to x^F}\frac{t^3k''(t)tk(t)}{(t^2k'(t))^2}
=\lim_{t\to x^F}\frac{(\theta^{-2}-3\theta^{-1}+2+o(1))(1+o(1))}{(\theta^{-1}-1+o(1))^2}=\frac{(\theta^{-1}-2)(\theta^{-1}-1)}{(\theta^{-1}-1)^2}
\end{array}
$$
and hence
$$
\begin{array}{c}
\dst\lim_{t\to x^F}\frac{\varphi'(t)}{k(t)(\varphi(t))^2}=\frac{1}{1-\theta}.
\end{array}
$$
Thus condition (\ref{gomes84}) holds for
 $\theta\not =1$ (the situation that we are always considering; see the last paragraph of section \ref{sint}). The rate of
convergence is given by $\gamma_n^2=(\gamma(\log
n))^2=(k'(b_n)/k^2(b_n))^2=((\theta-1)/\log n)^2$, i.e., also the
same order obtained in (\ref{gamatl}).
\end{rem}
%, also vanishing at $\theta=1$.\\

%Observe that the rate of convergence vanishes at $\theta=1$.\\

%Consider $u(t)=-\log(-\log F(t))$. Gomes (1989) imposes the
%condition
%$$
%\begin{array}{c}
%\dst\lim_{t\to x^F}\frac{\varphi'(t)}{u'(t)(\varphi(t))^2}=c<\infty
%\end{array}
%$$
%where $\varphi(t)=\nu''(t)/\nu'(t)$ with $\nu(t)=u^{-1}(t)$, and
%obtains boundedness of $(F^n(a_n
%x+b_n)-G_{\gamma_n}(x))/(\gamma_n-\gamma)^2$. Observe that
%$\nu'(t)=(u^{-1})'(t)=1/u'(u^{-1}(t))$ and
%$\nu''(t)=-u''(u^{-1}(t))(u^{-1})'(t)/(u'(u^{-1}(t)))^2$ as well
%as,
%$$
%\begin{array}{c}
%\dst\varphi(t)=\frac{\nu''(t)}{\nu'(t)}=-\frac{u''(u^{-1}(t))}{(u'(u^{-1}(t)))^2}
%\end{array}
%$$
%and
%$$
%\begin{array}{c}
%\dst\varphi'(t)=-\frac{u'''(u^{-1}(t))u'(u^{-1}(t))-2(u''(u^{-1}(t)))^2}{(u'(u^{-1}(t)))^4}.
%\end{array}
%$$
%Hence
%$$
%\begin{array}{c}
%\dst\frac{\varphi'(t)}{u'(t)(\varphi(t))^2}=\frac{-\frac{u'''(u^{-1}(t))u'(u^{-1}(t))-2(u''(u^{-1}(t)))^2}
%{(u'(u^{-1}(t)))^4}}{u'(t)\Big(\frac{u''(u^{-1}(t))}{(u'(u^{-1}(t)))^2}\Big)^2}=
%\frac{2}{u'(t)}-\frac{u'''(u^{-1}(t))u'(u^{-1}(t))}
%{u'(t)(u''(u^{-1}(t)))^2}
%\end{array}
%$$

\end{document}